\documentclass{article}
\usepackage{amsmath,amsthm,amscd,amssymb} 

 \newtheorem{theorem}{Theorem}[section]
 \newtheorem*{theoremA}{Theorem A}

 \newtheorem{proposition}[theorem]{Proposition}
 \newtheorem{lemma}[theorem]{Lemma}
 \newtheorem{corollary}[theorem]{Corollary}
 \newtheorem{remark}[theorem]{Remark}
 \newtheorem{definition}[theorem]{Definition}
 
 \newtheorem{example}[theorem]{Example}

 \numberwithin{equation}{section}
 \def\subrel#1#2{\mathrel{\mathop{#2}\limits_{#1}}}
 \def\qqq{\qed\vspace{2mm}} 

 \def\fp#1{{\mathbb F}_{\! p^{#1}}}

 \def\z2{{\mathbb Z}/2}

 \def\hz2{H\z2}

 \def\power#1{[\![#1]\!]}

 \def\im{\mbox{\rm Im}}

 \def\gal#1{\mbox{\rm Gal}(#1)}

 \def\inverselimit#1{\ \subrel{\subrel{#1}{\longleftarrow}}{\lim}}
 \def\directlimit#1{\ \subrel{\subrel{#1}{\longrightarrow}}{\lim}}

 \def\subrel#1#2{\mathrel{\mathop{#2}\limits_{#1}}}
 \def\qqq{\qed\vspace{2mm}} 
 \def\power#1{[\![#1]\!]}
\def\bigG{\mathbb G}

\begin{document}

\title{Comparison of Morava $E$-theories}
\author{Takeshi Torii
\footnote{
{\it Date}\,: October 14, 2008}
\setcounter{footnote}{-1}
\footnote{\textup{2000} {\it Mathematics Subject Classification}.
Primary 55N22; Secondary 55N20, 55S05.}
\setcounter{footnote}{-1}
\footnote{{\it Key words and phrases}.
Morava $E$-theory, formal group,
Morava stabilizer group, Chern character.}
}

\date{}
\maketitle

\begin{abstract}

In this note 
we show that 
the $n$th Morava $E$-cohomology group of a finite spectrum
with action of the $n$th Morava stabilizer group
can be recovered from 
the $(n+1)$st Morava $E$-cohomology group
with action of the $(n+1)$st 
Morava stabilizer group.
\end{abstract}

\section{Introduction}

This note is a continuation of 
\cite{Torii1}.
Let $p$ be a fixed prime number.
In the stable homotopy category ${\mathcal S}_{(p)}$ of $p$-local spectra,
there is a filtration of full subcategories ${\mathcal S}_n$,
where the objects of ${\mathcal S}_n$ consist of $E(n)$-local spectra.
The differences of each step of this filtration
are equivalent to the $K(n)$-local categories.
So it can be considered that 
the stable homotopy category ${\mathcal S}_{(p)}$
is built up from $K(n)$-local categories for various $n$.
In fact, the chromatic convergence theorem (cf. \cite{Ravenel-2})
says that the tower
$\cdots\to L_{n+1}X\to L_n X\to\cdots\to L_0X$ 
recovers a $p$-local finite spectrum $X$, that is,  
$X$ is homotopy equivalent to the
homotopy inverse limit of the tower.
Furthermore,
the chromatic splitting conjecture (cf. \cite{Hovey}) implies that
the $p$-completion of a finite spectrum $X$ is a direct summand of 
the product $\prod_n L_{K(n)}X$.
So it is important to know the relationship among
$K(n)$-local categories for various $n$.

We recall the classical Chern character map.
The Chern character map is a multiplicative natural transformation
from $K$-theory to rational cohomology:
\[ ch: K^*(X)\longrightarrow H^*(X;{\mathbb Q}[w^{\pm 1}])\cong
       \prod_{i\in {\mathbf Z}}H^{*+2i}(X;{\mathbb Q}). \]
The formal group law associated with $K$-theory
is the multiplicative one,
and the formal group law associated with the
rational cohomology is the additive one.
So the Chern character map is 
a multiplicative natural transformation of cohomology theories
from a height $1$ theory to a height $0$ theory.
Ando, Morava and Sadofsky~\cite{AMS}
have constructed a generalization of Chern character map
from the chromatic point of view.
In this note 
we construct a refinement of their 
generalization of Chern character map.
More precisely,
we construct a multiplicative 
natural transformation of cohomology theories
\[ \Theta: E_{n+1}^*(X)\longrightarrow \mathbb{B}^*(X). \] 
The cohomology theory $\mathbb{B}^*(X)$ is a coefficient extension 
of $E_n^*(X)$ so that
\[ \mathbb{B}^*(X)\cong \mathbb{B}^*\widehat{\otimes}_{E_n^*}
   E_n^*(X). \]
The natural transformation $\Theta$ can be regarded 
as a generalization of Chern character map
since
it is a multiplicative natural transformation of cohomology
theories from a height $(n+1)$ theory 
to a height $n$ theory.
There are natural actions of the Morava stabilizer groups
$G_n$ and $G_{n+1}$ on $\mathbb{B}^*(X)$,
which coincide on the Galois group.
Then the invariant submodule
$H^0(S_{n+1};\mathbb{B}^*(X))$ supports a natural structure
of $E_n^*$-module with compatible continuous $G_n$-action.
The following is the main theorem of this note.

\begin{theoremA}[Theorem~\ref{thm:invariantS_n+1}]\label{intro-main}
For every spectrum $X$,
there is a natural isomorphism of
$E_n^*$-modules with compatible continuous  
$G_n$-action
\[  E_n^*(X)\cong
    H^0(S_{n+1}; \mathbb{B}^*(X)).\]
For a finite spectrum $X$,
$E_n^*(X)$ with $G_n$-action can be recovered from
$E_{n+1}^*(X)$ with $G_{n+1}$-action.
\end{theoremA} 

The organization of this note is as follows:
In \S\ref{FGL}
we review the Lubin-Tate's deformation theory of formal group laws
and the results of \cite{Torii1}.
In \S\ref{stable-operations}
we review the relationship between the stable natural transformations of 
even-periodic complex oriented cohomology theories
and the homomorphisms of their formal group laws.
In \S\ref{lift-to-ch-0}
we construct a cohomology theory $\mathbb{B}^*(-)$
and the generalized Chern character map $\Theta$.
In \S~\ref{section:comparison:Galois-extension}
we show that the continuous $S_n$-cohomology 
with coefficients in $\mathbb{B}^*(X)$
is naturally isomorphic to $(L_{K(n)}E_{n+1})^*(X)$.
In \S\ref{section:comparison->En}
we prove the main theorem (Theorem~A) and some generalization.

\section{Formal group laws}\label{FGL}

In this section we review the deformation theory of
formal group laws.
In the following of this note a formal group law means 
a one-dimensional commutative formal group law.


Let $R_1$ and $R_2$ be two (topological) commutative rings.
Let $F_1$ (resp. $F_2$) be a 
formal group law over $R_1$ (resp. $R_2$).
We understand that a homomorphism from $(F_1,R_1)$ 
to $(F_2,R_2)$ 
is a pair $(f,\alpha)$ of a (continuous) ring homomorphism
$\alpha:R_2\to R_1$ and a homomorphism 
$f:F_1\to \alpha^*F_2$ in the usual sense,
where $\alpha^*F_2$ is the formal group law obtained from $F_2$
by the base change induced by $\alpha$.
We denote the set of all such pairs by 
\[ \mathbf{FGL}((F_1,R_1), (F_2,R_2)) .\]
If $R_1$ and $R_2$ are topological rings,
then we denote the subset of $\mathbf{FGL}((F_1,R_1),(F_2,R_2))$
consisting of $(f,\alpha)$ such that $\alpha$ is continuous by
\[ \mathbf{FGL}^c((F_1,R_1), (F_2,R_2)) .\]
The composition of two homomorphisms $(f,\alpha):(F_1,R_1)\to (F_2,R_2)$
and $(\beta,g):(F_2,R_2)\to (F_3,R_3)$ is defined as 
$(\alpha^*g\circ f,\alpha\circ\beta):(F_1,R_1)\to (F_3,R_3)$:
\[ F_1\stackrel{f}{\longrightarrow}\alpha^*F_2
      \stackrel{\alpha^*g}{\longrightarrow}
      \alpha^*(\beta^*F_3)=(\alpha\circ\beta)^*F_3.\]
A homomorphism $(f,\alpha):(F_1,R_1)\to (F_2,R_2)$ is an isomorphism
if there exists a homomorphism $(g,\beta):(F_2,R_2)\to (F_1,R_1)$ 
such that $(f,\alpha)\circ (g,\beta)=(X, id)$ and
$(g,\beta)\circ (f,\alpha)=(X, id)$.
Then a homomorphism $(f,\alpha):(F_1,R_1)\to (F_2,R_2)$ is an
isomorphism if and only if $\alpha$ is a (topological) ring isomorphism
and $f$ is an isomorphism in the usual sense. 

There is a $p$-typical formal group law $H_n$
over the prime field $\fp{}$ with $p$-series
\[ [p]^{H_n}(X)=X^{p^n},\]
which is called the height $n$ Honda
formal group law.
Let ${\mathbb F}$ be an algebraic extension of 
the finite field $\fp{n}$ with $p^n$ elements, 
and we suppose that  $H_n$ is defined over ${\mathbb F}$.
The automorphism group $S_n$ of $H_n$ over 
${\mathbb F}$ in the usual sense
is the $n$th Morava stabilizer group $S_n$, 
which is isomorphic to the unit group of
the maximal order of the central division algebra
over the $p$-adic number field ${\mathbf Q}_p$
with invariant $1/n$.
We denote by $G_n(\mathbb{F})$ the automorphism group of $H_n$
over ${\mathbb F}$ in the above sense:
\[ G_n(\mathbb{F})=\mathbf{Aut}(H_n, {\mathbb F}).\]
%
The automorphism group $G_n(\mathbb{F})$ is isomorphic to
the semi-direct product:
$G_n(\mathbb{F})\cong \Gamma(\mathbb{F})\ltimes S_n$,
where $\Gamma(\mathbb{F})$ is 
the Galois group $\gal{{\mathbb F}/\fp{}}$.

We recall Lubin and Tate's
deformation theory of formal group laws~\cite{Lubin-Tate}.
Let $R$ be a complete Noetherian local ring with 
maximal ideal $I$ such that
the residue field $k=R/I$ is of characteristic $p>0$.
Let $G$ be a formal group law over $k$ of height $n<\infty$.
Let $A$ be a complete Noetherian local $R$-algebra with maximal
ideal ${\mathfrak m}$. 
We denote by $\iota$ the canonical inclusion of residue fields
$k\subset A/{\mathfrak m}$ induced by the $R$-algebra structure.
A deformation of $G$ to $A$ is a 
formal group law $\widetilde{G}$ over $A$
such that $\iota^*G=\pi^*\widetilde{G}$
where $\pi:A\to A/{\mathfrak m}$ is the canonical projection.
Let $\widetilde{G}_1$ and $\widetilde{G}_2$ 
be two deformations of $G$ to $A$.
We define a $*$-isomorphism between $\widetilde{G}_1$ and 
$\widetilde{G}_2$ as an isomorphism
$\widetilde{u}:\widetilde{G}_1\to\widetilde{G}_2$ over $A$
such that $\pi^*\widetilde{u}$ is the identity map
between $\pi^*\widetilde{G}_1=\iota^*G=\pi^*\widetilde{G}_2$.
Then it is known that there is at most one $*$-isomorphism
between $\widetilde{G}_1$ and $\widetilde{G}_2$.
We denote by ${\mathcal C}(R)$ the category of complete Noetherian
local $R$-algebras with local $R$-algebra homomorphisms as morphisms. 
For an object $A$ of ${\mathcal C}(R)$,
we let $\mathbf{DEF}(A)$ be the set of all $*$-isomorphism classes of 
the deformations of $G$ to $A$.
Then $\mathbf{DEF}$ defines a functor from ${\mathcal C}(R)$
to the category of sets.
Let $R\power{t_i}=R\power{t_1,\ldots,t_{n-1}}$ be a formal power series ring 
over $R$ with $n-1$ indeterminates.
Note that $R\power{t_i}$ is an object of ${\mathcal C}(R)$.
Lubin and Tate constructed a formal group law 
$F(t_i)=F(t_1,\ldots,t_{n-1})$ over $R\power{t_i}$
such that 
for every deformation $\widetilde{G}$ of $G$ to $A$,
there is a unique local $R$-algebra homomorphism 
$\alpha:R\power{t_i}\to A$
such that $\alpha^*F(t_i)$
is $*$-isomorphic to $\widetilde{G}$. 
Hence 
the functor $\mathbf{DEF}$ is represented by
$R\power{t_i}$:
\[ \mathbf{DEF}(A)\cong \mbox{\rm Hom}_{{\mathcal C}(R)}
   (R\power{t_i},A)\]
and $F(t_i)$ is a universal object.


\begin{lemma}\label{key-lemma-ch0}
Let $F$ and $G$ be formal group laws of height $n<\infty$
over a field $k$ of characteristic $p>0$ and  
$(\overline{f},\overline{\alpha})$ 
an isomorphism from $(F,k)$ to $(G,k)$.
Let $R$ be a complete Noetherian local ring with residue field $k$
and $\alpha$ a ring automorphism of $R$ such that
$\alpha$ induces $\overline{\alpha}$ on the residue field.
Let $\widetilde{F}$ (resp. $\widetilde{G}$)
be a universal deformation of $F$ (resp. $G$)
over $R\power{u_i}=R\power{u_1,\ldots,u_{n-1}}$
(resp. $R\power{w_i}=R\power{w_1,\ldots,w_{n-1}}$).
Then there is a unique isomorphism $(g,\beta)$
from $(\widetilde{F},R\power{u_i})$ to $(\widetilde{G},R\power{w_i})$
such that 
$(g,\beta)$ induces $(\overline{f},\overline{\alpha})$
on the residue field
and $i\circ \alpha=\beta \circ j$,
where $i:R\to R\power{u_i}$ and $j:R\to R\power{w_i}$
are canonical inclusions. 
\end{lemma}

\proof
First, we show that there is such a homomorphism.
Let $f(X)\in R[X]$ be a lift of $\overline{f}(X)\in k[X]$
such that $f(0)=0$. 
Set $F'(X,Y)=f(\widetilde{F}(f^{-1}(X),f^{-1}(Y)))$.
Then $(F',R\power{u_i})$ is a deformation of
$\overline{\alpha}^*G$. 
We denote by $R'\power{u_i}$
the ring $R\power{u_i}$ with the $R$-algebra structure
given by $R\stackrel{\alpha}{\to}R\stackrel{i}{\to}R\power{u_i}$.
Then $(F', R'\power{u_i})$ is a deformation of $G$.
Since $\widetilde{G}$ is a universal deformation of $G$,
there exists a continuous $R$-algebra homomorphism 
$\beta: R\power{w_i}\to R'\power{u_i}$
and a $*$-isomorphism $\widetilde{u}: F'\to \beta^*\widetilde{G}$.
Then $(g,\beta)=(\widetilde{u}\circ f, \beta): 
(\widetilde{F}, R'\power{u}_i)\to (\widetilde{G}, R\power{w_i})$
is a lift of $(\overline{f},\overline{\alpha}): (F,k)\to (G,k)$.

By the same way,
we can construct a lift $(h,\gamma)$ of 
$(\overline{f},\overline{\alpha})^{-1}$.
Then $(h,\gamma)\circ (g,\beta)$ is a lift of 
$(X,id): (F,k)\to (F,k)$.
Note that $\beta\circ\gamma: R\power{u_i}\to R\power{u_i}$
is a continuous $R$-algebra homomorphism.   
Since $(\widetilde{F}, R\power{u_i})$ is a universal deformation,
$(h,\gamma)\circ (g,\beta)=(X,id)$ by the uniqueness.
Similarly,
we obtain that $(g,\beta)\circ (h,\gamma)=(X,id)$.
Hence we see that $(g,\beta)$ is an isomorphism
and a unique lift of 
$(\overline{f},\overline{\alpha}): (F,k)\to (G,k)$.
\qqq

Let ${\mathbf F}$ be an algebraic extension of $\fp{}$
which contains $\fp{n}$ and $\fp{n+1}$.
Let $W=W({\mathbf F})$ be the ring of Witt vectors
with coefficients in ${\mathbf F}$.
We define $E_n^0$ to be the ring of 
formal power series over $W$
with $(n-1)$ indeterminates:
\[ E_n^0=W\power{w_1,\ldots,w_{n-1}}.\]
The ring $E_n^0$ is a complete Noetherian local ring 
with residue field ${\mathbf F}$.
There is a $p$-typical formal group law $\widetilde{F}_n$
over $E_n^0$ with the $p$-series:
\begin{equation}\label{p-series-tilde-Fn+1}
 [p]^{\widetilde{F}_n}(X)=pX+_{\widetilde{F}_n}w_1X^p
                            +_{\widetilde{F}_n}w_2X^{p^2}
                            +_{\widetilde{F}_n}\cdots
                            +_{\widetilde{F}_n}w_{n-1}X^{p^{n-1}}
                            +_{\widetilde{F}_n}X^{p^n}.
\end{equation}
The formal group law $\widetilde{F}_n$ is a deformation of $H_n$
to $E_n^0$.  
Furthermore, 
$(\widetilde{F}_n, E_n^0)$ is a universal deformation 
of $(H_n,{\mathbf F})$.

\begin{lemma}\label{auto-gn}
The automorphism group $\mathbf{Aut}^c(\widetilde{F}_n,E_n^0)$ 
is isomorphic to $G_n(\mathbf{F})$.
\end{lemma}

We define $E_{n+1}^0$ to be 
a formal power series ring over $W$
with $n$ indeterminates:
\[ E_{n+1}^0=W\power{u_1,\ldots,u_n},\]
and there is a universal deformation
$(\widetilde{F}_{n+1}, E_{n+1}^0)$ of 
the height $(n+1)$ Honda group law $(H_{n+1},{\mathbf F})$. 
Let $R=E_{n+1}^0/I_n={\mathbf F}\power{u_n}$,
where $I_n=(p,u_1,\ldots,u_{n-1})$.
Let $F_{n+1}=\pi^*\widetilde{F}_{n+1}$,
where $\pi$ is the quotient map $E_{n+1}^0\to R$.
Then $F_{n+1}$ is a deformation of $H_{n+1}$ to $R$.
The following lemma is easy.

\begin{lemma}\label{automorphism_group_F_n+1}
The automorphism group $\mathbf{Aut}^c(F_{n+1},R)$
is isomorphic to $G_{n+1}(\mathbf{F})$.
\end{lemma}

If we suppose that $F_{n+1}$ is defined 
over the quotient field ${\mathbf F}((u_n))$ of $R$,
then its height is $n$.
Since the formal group laws over a separably closed field
is classified by their height,
there is an isomorphism $\Phi$ between $F_{n+1}$ and $H_n$
over the separable closure ${\mathbf F}((u_n))^{sep}$ of 
${\mathbf F}((u_n))$ 
(cf. \cite{Lazard, Hazewinkel}).
We fix such an isomorphism $\Phi$.
Since $\Phi: F_{n+1}\to H_n$ is a homomorphism between $p$-typical formal
group laws, 
$\Phi$ has a following form:
\[ \Phi(X)=\sum_{i\ge 0}{}^{H_n}\Phi_i X^{p^i}.\]
Let $L$ be the extension field of ${\mathbf F}((u_n))$
obtained by adoining all the coefficients of the isomorphism $\Phi$.
So $(\Phi,id_L)$ is
an isomorphism from $(F_{n+1},L)$ to $(H_n,L)$:
\[   (\Phi, id_L): (F_{n+1},L)\stackrel{\cong}{\longrightarrow}
                   (H_n,L).\]
Note that $L$ is a totally ramified Gaois extension
of infinite degree over ${\mathbf F}((u_n))$ 
with Galois group isomorphic to $S_n$
(\cite{Gross, Torii1}).

In the following of this note 
we abbreviate 
$G_n=G_n(\mathbf{F})$, $G_{n+1}=G_{n+1}(\mathbf{F})$
and $\Gamma=\Gamma(\mathbf{F})$.
There are quotient maps 
$G_n\to\Gamma$
and $G_{n+1}\to\Gamma$.
We define $\mathbb{G}$ to be the fibre product
of $G_n$ and $G_{n+1}$ over $\Gamma$:
$\mathbb{G}=G_n\times_{\Gamma}G_{n+1}$.
Note that there is an isomorphism 
$\mathbb{G}\cong \Gamma\ltimes (S_n\times S_{n+1})$.
Then $G_{n+1}\cong
\Gamma\ltimes S_{n+1}$ and $G_n=\Gamma\ltimes S_n$
are subgroups of ${\mathbb G}$.
In \cite{Torii1} we have shown the following theorem.

%
%
%

\begin{theorem}[{cf. \cite[\S2.4]{Torii1}}]\label{FGL-main-theorem-v3}
The pro-finite group ${\mathbb G}$ acts on 
$(F_{n+1},L)\stackrel{}{\cong} (H_n,L)$.
The action of the subgroup $G_{n+1}$ on $(F_{n+1},L)$
is an extension of the action on $(F_{n+1},R)$,
and the action of the subgroup $G_n$ on $(H_n,L)$
is an extension of the action on $(H_n,{\mathbf F})$.
\end{theorem}

\section{Natural transformations of cohomology theories}
\label{stable-operations}

Let $\mathcal{S}$ be the stable homotopy category.
For a spectrum $h\in\mathcal{S}$,
we denote by $h^*(-)$ the associated 
generalized cohomology theory.
For spectra $h$ and $k$,
we denote by ${\mathbf C}(h,k)$ 
the set of all degree $0$
natural transformation of cohomology theories from $h$ to $k$.
Then ${\mathbf C}(h,k)$ is naturally identified with
the set of all degree $0$ morphisms from $h$ to $k$
in $\mathcal{S}$.

\begin{definition}\rm
Let $R$ be a 
commutative ring.
A topological $R$-module $M$ is said to be 
{\it linearly topologized} if $M$ has a fundamental 
neighbourhood system at the zero consisting 
of the open submodules.
A linearly topologized $R$-module $M$ 
is said to be {\it linearly compact}
if it is Hausdorff and it has the finite intersection
property with respect to the closed cosets
A topological ring $R$ is linearly compact
if $R$ is linearly compact as an $R$-module.
(cf. 
\cite[Definition~2.3.13]{HPS}).
\end{definition}

\begin{example}\rm
A linearly topologized compact Housdorff (e.g. profinite) 
module is linearly compact.
If $R$ is a complete Noetherian local ring,
then a finitely generated $R$-module is linearly compact.
In particular,
a finite dimensional vector space over a field 
is linearly compact. 
\end{example}

\begin{lemma}[cf. {\cite[Corollary~2.3.15]{HPS}}]
\label{exact-inverse-limit}
Let ${\mathcal I}$ be a filtered category.
The inverse limit functor indexed by ${\mathcal I}$
is exact
in the category of linearly compact modules
and continuous homomorphisms.
\end{lemma}

For a graded commutative ring $h^*$,
we say that $h^*$ is even-periodic 
if there is a unit $u\in h^{-2}$ of degree $(-2)$ and $h^{\rm odd}=0$.
Note that $h^*=h^0[u^{\pm 1}]$ if $h^*$ is even-periodic.
We say that a commutative ring spectrum $h$ is even-periodic
if the coefficient ring $h^*=h^*(S^0)$ is even-periodic.
For a spectrum $X\in\mathcal{S}$,
let $\Lambda(X)$ be the full subcategory of the comma category
$(\mathcal{S}\downarrow X)$,
whose objects are maps $X_{\lambda}\to X$ 
with $X_{\lambda}$ finite.
Then $\Lambda(X)$ is an essentially small filtered category 
(see \cite[Definition~2.3.7]{HPS}).

\begin{lemma}\label{lemma:inverselimit-description}
Suppose that $h$ is an even-periodic commutative ring spectrum
with $h^0$ Noetherian and linearly compact.
Then there is a natural isomorphism
\[  h^*(X)\stackrel{\cong}{\longrightarrow}\
    \subrel{\longleftarrow}{\lim}\ h^*(X_{\lambda})\]
for any spectrum $X$,
where the inverse limit on the right hand side is taken 
over $\Lambda(X)$.
\end{lemma}

\proof
Since 
$h^*(X_{\lambda})$ is a finitely generated module over $h^*$
for a finite spectrum $X_{\lambda}$,
it is a linearly compact module.
Then 
${\lim}\,h^*(X_{\lambda})$
is a cohomology theory by Lemma~\ref{exact-inverse-limit}
(cf. \cite[Proposition~2.3.16]{HPS}).
The natural transformation 
$h^*(X)\to{\lim}\,h^*(X_{\lambda})$
of cohomology theories gives
an isomorphism for $X=S^0$.
Therefore it is an isomorphism for all $X$.
\qqq

By Lemma~\ref{lemma:inverselimit-description},
if $h$ is even-periodic and $h^0$ is Noetherian 
linearly compact,
$h^*(X)$ naturally supports a structure of 
linearly compact module over $h^*$
for all $X$.

Let $p$ be a prime number and 
$BP$ the Brown-Peterson spectrum at $p$.
Let $h^*$ be a graded commutative ring over ${\mathbb Z}_{(p)}$.
We suppose that there is a $p$-typical 
formal group law $F_h'$ of degree $(-2)$ over $h^*$.
Since the associated formal group law $F_{BP}'$ to $BP$ 
is universal with respect to $p$-typical ones,
there is a unique ring homomorphism $r: BP^*\to h^*$.
Then the functor $h^*\otimes_{BP^*}BP^*(-)$
is a generalized cohomology theory on the category of finite spectra
if $p,v_1,v_2,\ldots$ is a regular sequence in $h^*$
by the Landweber exact functor theorem~\cite{Landweber}. 
We say such a graded commutative ring $h^*$ is
Landweber exact over $BP^*$.
For any spetrum $X$, Lemma~\ref{lemma:inverselimit-description} 
suggests to define 
\[ \begin{array}{rcl}
     h^*(X)&=&h^*\widehat{\otimes}_{BP^*}BP^*(X)\\[2mm]
           &=&\inverselimit{}\left(h^*
              \otimes_{BP^*}BP^*(X_{\lambda})\right), 
   \end{array}\] 
where the inverse limit is taken over $\Lambda(X)$.

\begin{lemma}\label{cohomology}
Suppose that $h^*$ is an even-periodic
Landweber exact graded commutative ring over $BP^*$.
Furthermore, suppose that $h^0$ is Noetherian and linearly compact.
Then the functor $h^*(-)$
is a complex oriented commutative multiplicative cohomology theory.
\end{lemma}

\proof
It is easy to see that $h^*(-)$ takes coproducts to products.
The exactness of $h^*(-)$ follows from 
Lemma~\ref{exact-inverse-limit}
(cf. \cite[Proposition~2.3.16]{HPS}).
Hence $h^*(-)$ is a commutative multiplicative cohomology theory.
The natural transformation $BP^*(-)\to h^*(-)$
gives us an orientation of $h^*(-)$.
\qqq

\begin{definition}\rm
Let $h$ and $k$ be commutative ring spectra.
We denote by $\mathbf{Mult}(h,k)$ the set of
all multiplicative natural transformations of cohomology
theories from $h^*(-)$ to $k^*(-)$.
If $h^*(-)$ and $k^*(-)$ have their values in the category 
of linearly compact modules,
then we denote by $\mathbf{Mult}^c(h,k)$ 
the subset of $\mathbf{Mult}(h,k)$ consisting
of $\theta$ such that $\theta: h^*(X)\to k^*(X)$ 
is continuous for all $X$.
\end{definition}

If $h^*(-)$ is a complex oriented cohomology theory,
then the orientation class $X_h\in h^2({\mathbf C}P^{\infty})$
gives a formal group law $F_h'$ of degree $2$.
Furthermore, if $h$ is even-priodic,
then a unit $u\in h^{-2}$ gives a degree $0$ formal group law
by $F_h(X,Y)= u F_h'(u^{-1}X, u^{-1}Y)$.
In the following of this section
we suppose that $h$ and $k$ are even-periodic
complex orientable commutative ring spectra. 
Furthermore, we fix a unit $u\in h^{-2}$ (resp. $v\in k^{-2}$)
and an orientation class $X_h\in h^2({\mathbf C}P^{\infty})$
(resp. $X_k\in k^2({\mathbf C}P^{\infty})$). 
Then we obtain a degree $0$ formal group law
$F_h$ (resp. $F_k$) associated with $h$ (resp. $k$)
as above.
A multiplicative natural transformation $\theta: h^*(-)\to k^*(-)$
gives a ring homomorphism $\alpha: h^0\to k^0$ and
an isomorphism $f: F_k\to \alpha^*F_h$ of 
formal group laws.
Note that $f(X)=\theta(u)\widetilde{f}(v^{-1}X)$,
where $\widetilde{f}(X_k)=\theta(X_h)$.
In particular,
$f'(0)=\theta (u) v^{-1}$ is a unit of $k^0$.
Hence we obtain a map from $\mathbf{Mult}(h,k)$
to $\mathbf{FGL}((F_k,k^0),(F_h,h^0))$:
\[ \Xi:  \mathbf{Mult}(h,k) \to \mathbf{FGL}((F_k,k^0),(F_h,h^0)) .\] 
If $h^0$ and $k^0$ are Noetherian and linearly compact,
then $\Xi$ induces 
\[ \Xi^c: \mathbf{Mult}^c(h,k) \to \mathbf{FGL}^c((F_k.k^0),(F_h,h^0)) .\]   

\begin{remark}\label{cotinuity-remark}\rm
Let $h^0$ and $k^0$ be Noetherian linearly compact rings. 
If $\theta\in\mathbf{Mult}(h,k)$ induces a continuous
ring homomorphism $h^0\to k^0$,
then $\theta\in\mathbf{Mult}^c(h,k)$.
\end{remark}

\begin{proposition}\label{operation-FGL}
Suppose that $h^*$ and $k^*$ are even-periodic
Landweber exact graded commutative rings over $BP_*$. 
Then the map
$\Xi:\mathbf{Mult}(h,k)\longrightarrow
\mathbf{FGL}((F_k,k_0),(F_h,h_0))$
is a bijection.
Furthermore,
if $h^0$ and $k^0$ are Noetherian and linearly compact,
then $\Xi^c:\mathbf{Mult}^c(h,k)\longrightarrow
\mathbf{FGL}^c((F_k,k_0),(F_h,h_0))$
is also a bijection.
\end{proposition}

\proof
This follows from the Landweber exact functor theorem
and the fact that the graded commutative ring $BP_*(BP)$ 
represents the set of all triples
$(F,f,G)$ where $F,G$ are $p$-typical formal group laws
and $f$ is a strict isomorphism betweent them. 
\qqq

\section{The generalized Chern character $\Theta$}\label{lift-to-ch-0}

In this section we construct 
the generalized Chern character $\Theta$,
which is a multiplicative natural transformation
of cohomology theories
from the height $(n+1)$ cohomology theory $E_{n+1}^*(-)$
to the height $n$ cohomology theory $\mathbb{C}^*(-)$,
which is a coefficient extension of
$E_n^*(-)$.

We recall that ${\mathbf F}$ is an algebraic extension
of $\fp{}$ which contains the finite fields $\fp{n}$ and $\fp{n+1}$.
Also recall that
$E_n^0=W\power{w_1,\ldots,w_{n-1}}$
and 
$E_{n+1}^0=W\power{u_1,\ldots,u_n}$,
where $W$ is the ring of Witt vectors with coefficients in $\mathbf{F}$.
We define graded commutative rings $E_n^*$ and $E_{n+1}^*$ by
$E_n^*=E_n^0[u^{\pm 1}]$ and $E_{n+1}^*=E_{n+1}^0[u^{\pm 1}]$,
where $|u|=|w|=-2$.
Hence we have 
\[ \begin{array}{rcl}
    E_n^*&=  &W\power{w_1,\ldots,w_{n-1}}[w^{\pm 1}],\\[2mm]
    E_{n+1}^*&=&W\power{u_1,\ldots,u_n}[u^{\pm 1}].\\
   \end{array}\]
The $p$-typical formal group law $\widetilde{F}_n$
over $E_n^0$ gives a ring homomorphism
$r_n: BP^*\to E_n^*$, which is 
given by $r_n(v_i)=w_iw^{(p^i-1)}\ (1\le i<n), 
r_n(v_n)=w^{(p^n-1)}, r_n(v_i)=0\ (i>n)$.
Also, the $p$-typical formal group law $\widetilde{F}_{n+1}$
over $E_{n+1}^0$ give a 
ring homomorphism $r_{n+1}: BP^*\to E_{n+1}^*$,
which is given by
$r_{n+1}(v_i)=u_i u^{(p^i-1)}\ (1\le i\le n),
 r_{n+1}(v_{n+1})=u^{(p^{n+1}-1)}, 
 r_{n+1}(v_i)=0\ (i>n+1)$.
Then we can regard $E_n^*$ and $E_{n+1}^*$
as even-periodic Landweber exact graded commutative
rings over $BP^*$
through the ring homomorphisms $r_n$ and $r_{n+1}$, respectively.
Furthermore, since $E_n^0$ and $E_{n+1}^0$
are Noetherian local rings,
they are linearly compact.
Hence $E_n^*(-)=E_n^*\widehat{\otimes}_{BP^*}BP^*(-)$ and 
$E_{n+1}^*(-)=E_{n+1}^*\widehat{\otimes}_{BP^*}BP^*(-)$
are complex oriented commutative multiplicative 
cohomology theories by Lemma~\ref{cohomology}.
Let $A=(W((u_n)))^{\wedge}_p$ be the $p$-adic completion of 
$W((u_n))$, the ring of Laurent series over $W$. 
Then $A$ is a complete discrete valuation ring
with uniformizer $p$   
and residue field $M={\mathbf F}((u_n))$. 
In particular,
$A$ is a Henselian ring. 
We recall the following lemma on Henselian rings.

\begin{lemma}[{cf.~\cite[Proposition~I.4.4.]{Milne}}]\label{henselian}
Let $R$ be a Henselian ring with residue field $k$.
Then the functor $R\mapsto S\otimes_R k$ induces
an equivalence between the category of finite \'{e}tale  
$R$-algebras and the category of finite \'{e}tale $k$-algebras.
\end{lemma}

In \cite[\S2.3]{Torii1}
we have constructed a sequence of finite Galois extensions of 
${\mathbf F}((u_n))$:
\begin{align}\label{sequence-Galois-extesion-residue} 
   {\mathbf F}((u_n))=L({-1})\to L(0)\to L(1)\to\cdots,
\end{align}
where $L(i)$ is obtained by adjoining the coefficients 
$\Phi_0,\Phi_1,\ldots,\Phi_i$ of the isomorphism 
$\Phi: F_{n+1}\stackrel{\cong}{\to} H_n$
of formal group laws. 
We denote by $S_n(i)$
the Galois group for $L(-1)\to L(i)$.
Then $S_n(i)$ is a finite group of order $(p-1)p^{i}$,
and $S_n=\ \subrel{\longleftarrow i}{\lim}\, S_n(i)$.
We let $S_n^{(i)}$ be the kernel of the canonical 
surjection $S_n\to S_n(i)$.
By definition,
$L=\ \subrel{\longrightarrow i}{\rm lim}L(i)
=\cup_i L(i)$, and
${\mathbf F}((u_n))\to L$ is an infinite Galois extension
with Galois group $S_n$.
Note that $L(i)$ is stable under the action of $\bigG$
for all $i$.

By Lemma~\ref{henselian},
we obtain a sequence of finite \'{e}tale 
$A$-algebras:
\[ A=B({-1})\to B(0)\to B(1)\to \cdots .\] 
We denote by $B({\infty})$ the direct limit 
$\subrel{\longrightarrow i}{\rm lim}B(i)$
and $B$ the $p$-adic completion of $B({\infty})$.

\begin{lemma}       
The ring $B$ is a complete discrete valuation ring 
of characteristic $0$ 
with uniformizer $p$ and residue field 
$L=\ \subrel{\longrightarrow i}{\rm lim}L(i)$.
\end{lemma}

\proof
Since $L(i)$ is a separable extension over ${\mathbf F}((u_n))$,
we can take $a\in L(i)$ such that $L(i)={\mathbf F}((u_n))(a)$.
Let $f(X)\in {\mathbf F}((u_n))[X]$ be the minimal polynomial of $a$
and $\widetilde{f}(X)\in A[X]$ a monic polynomial
which is a lift of $f(X)$.
Then $B(i)\cong A[X]/(\widetilde{f}(X))$.
Then we see that $B(i)$ is a complete 
discrete valuation ring with uniformizer $p$
and residue field $L(i)$.
This implies that 
$B({\infty})$ is also a discrete valuation ring with uniformizer $p$
and residue field $L$.
Then the lemma follows from the fact that
$B$ is the $p$-adic completion of $B({\infty})$.
\qqq

We abbreviate $B\power{w_1,\ldots,w_{n-1}}$
and $B\power{u_1,\ldots,u_{n-1}}$ by
$B\power{w_i}$ and $B\power{u_i}$, respectively, etc.
Then we obtain a sequence of finite \'{e}tale $A\power{u_i}$-algebras:
\[ A\power{u_i}= B({-1})\power{u_i}
               \to B(0)\power{u_i}
               \to B(1)\power{u_i}\cdots, \]
and 
$B\power{u_i}$ is the $I_n$-adic
completion of 
$\subrel{\longrightarrow j}{\rm lim}B(j)\power{u_i}$,
where $I_n=(p,u_1,\ldots,u_{n-1})$.

We define an even-periodic graded commutative ring $\mathbb{A}^*$ by
\[ \mathbb{A}^*=A\power{u_1,\ldots,u_n}[u^{\pm 1}].\]
There is a canonical inclusion
$E_{n+1}^*\hookrightarrow \mathbb{A}^*$, and 
the ring homomorphism
$BP^*\to E_{n+1}^*\to \mathbb{A}^*$
satisfies the Landweber exact condition.
Also, the degree $0$ subring $\mathbb{A}^0=A\power{u_i}$
is a complete Noetherian local ring.
Hence $\mathbb{A}^*(-)=
\mathbb{A}^*\widehat{\otimes}_{BP^*}BP^*(-)$
gives a complex oriented commutative multiplicative cohomology theory
by Lemma~\ref{cohomology}.
We denote by $\mathbb{A}$ the representing commutative rings spectrum
of $\mathbb{A}^*(-)$.
Then the canonical inclusion $E_{n+1}^*\hookrightarrow \mathbb{A}^*$
induces a ring spectrum map $E_{n+1}\to\mathbb{A}$.

\begin{lemma}\label{lemma:A=L_{K(n)}E_{n+1}}
$\mathbb{A}$ is equivalent to $L_{K(n)}E_{n+1}$ 
as a commutative ring spectrum.
\end{lemma}

\proof
There is a tower $\{M(I)\}_I$ of generalized Moore spectra 
of height $n$
such that $L_{K(n)}X\simeq\ {\rm holim}_I\, L_nX\wedge M(I)$
for any spectrum $X$ (cf. \cite[Proposition~7.10(e)]{Hovey-Strickland}).
In paricular,
$L_{K(n)}E_{n+1}\simeq\ {\rm holim}_I\,L_nE_{n+1}\wedge M(I)$.
Then $\pi_*L_{K(n)}E_{n+1}\cong A\power{u_i}[u^{\pm 1}]$.
Hence we can identify the localization map
$E_{n+1}\to L_{K(n)}E_{n+1}$ with the ring spectrum map
$E_{n+1}\to\mathbb{A}$.
\qqq

We define an even-periodic graded commutative ring $\mathbb{B}^*$ by
\[ \mathbb{B}^*=B\power{u_1,\ldots,u_n}[u^{\pm 1}].\]
There is a canonical inclusion
$\mathbb{A}^*\hookrightarrow \mathbb{C}^*$, and 
the ring homomorphism
$BP^*\to \mathbb{A}^*\to \mathbb{B}^*$
satisfies the Landweber exact condition.
Also, the degree $0$ subring $\mathbb{B}^0=B\power{u_i}$
is a complete Noetherian local ring.
Hence $\mathbb{B}^*(-)=
\mathbb{B}^*\widehat{\otimes}_{BP^*}BP^*(-)$
gives a complex oriented commutative multiplicative cohomology theory
by Lemma~\ref{cohomology}.
We denote by $\mathbb{B}$ the representing commutative rings spectrum
of $\mathbb{B}^*(-)$.
Then the following lemma is clear from the construction.

\begin{lemma}
There is a natural isomorphism
$\mathbb{B}^*(X)\cong \mathbb{B}^*
\widehat{\otimes}_{\mathbb{A}^*}\mathbb{A}^*(X)
=\ \subrel{\longleftarrow\lambda}{\lim} 
\mathbb{B}^*\otimes_{\mathbb{A}^*}\mathbb{A}^*(X_{\lambda})$
for all spectra $X$,where the inverse limit is taken over
$\Lambda(X)$.
\end{lemma}

We define an even-periodic graded commutative ring $\mathbb{C}^*$ by
\[ \mathbb{C}^*=B\power{w_1,\ldots,w_{n-1}}[w^{\pm 1}].\]
There is a canonical inclusion
$E_n^*\hookrightarrow \mathbb{C}^*$, and 
the ring homomorphism
$BP^*\to E_n^*\to \mathbb{C}^*$
satisfies the Landweber exact condition.
Also, the degree $0$ subring $\mathbb{C}^0=B\power{w_i}$
is a complete Noetherian local ring.
Hence $\mathbb{C}^*(-)=
\mathbb{C}^*\widehat{\otimes}_{BP^*}BP^*(-)$
gives a complex oriented commutative multiplicative cohomology theory
by Lemma~\ref{cohomology}.
Then the following lemma is clear from the construction.

\begin{lemma}\label{lemma:description-C-theory}
There is a natural isomorphism
$\mathbb{C}^*(X)\cong \mathbb{C}^*\widehat{\otimes}_{E_n^*}E_n^*(X)
=\ \subrel{\longleftarrow\lambda}{\lim} 
\mathbb{C}^*\otimes_{E_n^*}E_n^*(X_{\lambda})$
for all spectra $X$, where the inverse limit is taken over
$\Lambda(X)$.
\end{lemma}

The ring homomorphism $E_{n+1}^*\to \mathbb{A}^*\to\mathbb{B}^*$
give a formal group law
$\widetilde{F}_{n+1}$ over $\mathbb{B}^0=B\power{u_i}$, and
the ring homomorphism $E_n^*\to \mathbb{C}^*$ gives
a formal group law $\widetilde{F_n}$ over $\mathbb{C}^0=B\power{w_i}$.

\begin{lemma}\label{universality}
The formal group laws $(\widetilde{F}_{n+1},B\power{u_i})$
and $(\widetilde{F}_n,B\power{w_i})$ are
universal deformations of
$(F_{n+1},L)$ and $(H_n,L)$,respectively,
on the category of complete Noetherian local $B$-algebras.
\end{lemma}

\proof
From the fact that $(\widetilde{F}_n,E_n^0)$ is a universal
deformation of $(H_n,{\mathbf F})$,
it is easy to see that 
$(\widetilde{F}_n, B\power{w_i})$ is a universal deformation of
$(H_n,L)$.
From the form of the $p$-series of $\widetilde{F}_{n+1}$ given by 
(\ref{p-series-tilde-Fn+1}),
we see that $(\widetilde{F}_{n+1},B\power{u_i})$
is a universal deformation of $(F_{n+1},L)$.
\qqq


\begin{corollary}
The action of $G_{n+1}$ on $(\widetilde{F}_{n+1},E_{n+1}^0)$
extends to an action on $(\widetilde{F}_{n+1},B\power{u_i})$
such that the induced action on $(F_{n+1},L)$ 
coincides with the action of Theorem~\ref{FGL-main-theorem-v3}.
\end{corollary}

\proof
It is sufficient to show that 
the action of $G_{n+1}$ on $E_{n+1}^0$ extends
to an action on $B\power{u_i}$.
For $g\in G_{n+1}$,
$u_n^g$ is a unit multiple of $u_n$ 
modulo $(p,u_1,\ldots,u_{n-1},u_n^2)$.
Hence the ring homomorphism
$E_{n+1}^0\stackrel{g}{\to}E_{n+1}^0\to
 (E_{n+1}^0[u_n^{-1}])^{\wedge}_{I_n}=A\power{u_i}$
extends to a ring homomorphism
$E_{n+1}[u_n^{-1}]\to A\power{u_i}$
This induces a ring homomorphism 
$A\power{u_i}\to A\power{u_i}$
and defines an action of $G_{n+1}$ on $A\power{u_i}$.
Since $B(j)\power{u_i}\to B({j+1})\power{u_i}$
is \'{e}tale for $j\ge -1$
and $L(j)$ is stable under the action of $G_{n+1}$ on $L$,
the action on $A\power{u_i}$ extends to 
$B(j)\power{u_i}$ uniquely and compatibly by Lemma~\ref{henselian}.
Hence we obtain an action on 
$\subrel{\longrightarrow j}{\rm lim}B(j)\power{u_i}$
and its $I_n$-adic completion $B\power{u_i}$.
\qqq


We denote the action of $G_{n+1}$ on 
$(\widetilde{F}_{n+1},B\power{u_i})$ 
by $\Upsilon(g)=(t(g), \upsilon(g)): 
(\widetilde{F}_{n+1},B\power{u_i})\to 
(\widetilde{F}_{n+1},B\power{u_i})$ for $g\in G_{n+1}$.

\begin{corollary}
The $(n+1)$th extended Morava stabilizer group 
$G_{n+1}$ acts on the cohomology theory
$\mathbb{B}^*(-)$ 
as multiplicative cohomology operations.
\end{corollary}

\proof
This follows from Proposition~\ref{operation-FGL}. 
\qqq

There are isomorphisms
$S_n\cong \gal{L/{\mathbf F}((u_n))}$ and 
$G_n\cong \gal{L/\fp{}((u_n))}$
through the action of $G_n$ on $L$
(cf. \cite{Gross, Torii1}). 

\begin{lemma}
The action of $G_n$ on $L$ lifts to the action on $B$.
\end{lemma}

\proof
Since $L(i)$ is stable under the action of $G_n$ on $L$
for all $i\ge -1$,
the action of $G_n$ on $L(i)$ lifts to the action on
$B(i)$ compatibly by Lemma~\ref{henselian}.
This induces an action on $B({\infty})$.
Since $B$ is the $p$-adic completion of $B({\infty})$,
we obtain an action on $B$
which is a lift of the action on $L$.
\qqq

We denote this action of $G_n$ on $B$ 
by $\tau(g): B\to B$
for $g\in G_n$.  
Since the actions of $G_n$ on $E_n^0$ and $B$
are compatible on $W$,
the diagonal action defines an action of $G_n$ on $B\power{w_i}=
B\widehat{\otimes}_W E_n^0$. 
Then we obtain an extension of the action of $G_n$ 
on $(\widetilde{F}_n, E_n^0)$ to $(\widetilde{F}_n, B\power{w_i})$.
We denote this action of $G_n$ on $(\widetilde{F}_n, B\power{w_i})$
by $\Omega(g)=(s(g),\omega(g)):
(\widetilde{F}_n, B\power{w_i})\to (\widetilde{F}_n, B\power{w_i})$
for $g\in G_n$.

\begin{corollary}
The $n$th extended Morava stabilizer group
$G_n$ acts on 
$\mathbb{C}^*(-)$ as 
multiplicative cohomology operations. 
\end{corollary}

\proof
This follows from Proposition~\ref{operation-FGL}. 
\qqq

\begin{lemma}
There is a unique isomorphism 
$(\widetilde{\Phi},\widetilde{\varphi}):
(\widetilde{F}_{n+1},B\power{u_i})\stackrel{\cong}{\to}
(\widetilde{F}_n,B\power{w_i})$
such that $\widetilde{\varphi}$ is a continuous $B$-algebra
homomorphism and $\widetilde{\Phi}$ induces
$\Phi$ on the residue fields.
\end{lemma}

\proof
Since there is an isomorphism
$(\Phi,id_L):(\widetilde{F}_{n+1},L)\to (H_n, L)$, 
the lemma follows from Lemma~\ref{key-lemma-ch0}.
\qqq

%

\begin{lemma}\label{fundamental-commutativity}
For $g\in G_n$,
there is a commutative diagram:
\[ \begin{array}{ccc}
    (\widetilde{F}_{n+1},B\power{u_i}) & 
    \stackrel{(X,\theta(g))}{\longrightarrow} &
    (\widetilde{F}_{n+1},B\power{u_i}) \\
    {\scriptstyle (\widetilde{\Phi},\widetilde{\varphi})}
    \bigg\downarrow & & \bigg\downarrow 
    {\scriptstyle (\widetilde{\Phi},\widetilde{\varphi})}\\
    (\widetilde{F}_n,B\power{w_i}) & 
    \stackrel{\Omega(g)}{\longrightarrow} &
    (\widetilde{F}_n,B\power{w_i}), \\
   \end{array}\]
where $\theta(g):B\power{u_i}\to B\power{u_i}$
is given by $\theta(g)(b)=\tau(g)(b)$ for $b\in B$
and $\theta(g)(u_i)=u_i$ for $1\le u_i<n$.
\end{lemma}

\proof
Note that $(\theta(g)\circ\widetilde{\varphi})|_B=
\tau(g)=(\widetilde{\varphi}\circ \omega(g))|_B$.
The diagram induced on the residue field is commutative
by definition of the action of $G_n$ on $(F_{n+1},L)\cong (H_n,L)$.
Then the lemma follows from the universality of
$(\widetilde{F}_n,B\power{w_i})$.
\qqq

\begin{lemma}\label{lemma:Gn-module-inclusion}
For $g\in G_{n+1}$,
there is a commutative diagram:
\[ \begin{array}{ccc}
    (\widetilde{F}_{n+1}, B\power{u_i}) & 
    \stackrel{\Upsilon(g)}{\longrightarrow} & 
    (\widetilde{F}_{n+1}, B\power{u_i}) \\
    {\scriptstyle (\widetilde{\Phi},\widetilde{\varphi})}
    \bigg\downarrow & & \bigg\downarrow 
    {\scriptstyle (\widetilde{\Phi},\widetilde{\varphi})}\\
    (\widetilde{F}_n, B\power{w_i}) & 
    \stackrel{(X,\mu(g))}{\longrightarrow} &
    (\widetilde{F}_n, B\power{w_i}),\\
   \end{array}\] 
where $\mu(g)$ is given by 
$\mu(g)(b)=\widetilde{\varphi}^{-1}(\upsilon(g)(b))$ for $b\in B$
and $\mu(g)(w_i)=w_i$ for $1\le i<n$. 
\end{lemma}

\proof
Note that $(v(g)\circ\widetilde{\varphi})|_B=
v(g)|_B=(\widetilde{\varphi}\circ \mu(g))|_B$.
The diagram induced on the residue field is commutative
by definition of the action of $G_{n+1}$ on $(F_{n+1},L)\cong (H_n,L)$.
Then the lemma follows from the universality of
$(\widetilde{F}_n,B\power{w_i})$.
\qqq

\begin{corollary}\label{cor:commutativity-actions}
The profinite group $\bigG$ acts on 
$(\widetilde{F}_{n+1},B\power{u_i})
\cong (\widetilde{F}_n,B\power{w_i})$ such that
the action of the subgroup $G_{n+1}$ coincides with
$\Upsilon$,
and the action of the subgroup $G_n$ coincides with 
$\Omega$.
\end{corollary}

\proof
We have the action $\Upsilon$ of $G_{n+1}$ on 
$(\widetilde{F}_{n+1},B\power{u_i})$ and 
the action $\Omega$ of $G_n$ on $(\widetilde{F}_n,B\power{w_i})$.
The action of the subgroup $\Gamma$ of $G_{n+1}$ 
on $(\widetilde{F}_{n+1},B\power{u_i})$ coincides with
the action on $(\widetilde{F}_n,B\power{w_i})$
as the subgroup of $G_n$ under the isomorphism
$(\widetilde{\Phi},\widetilde{\varphi})$.
Hence it is sufficient to show that 
the following diagram commutes for
$g\in S_{n+1}$ and $h\in S_n$:
\[  \begin{array}{ccc}
     (\widetilde{F}_{n+1},B\power{u_i}) &
     \stackrel{(X,\theta(h))}{\longrightarrow} &
     (\widetilde{F}_{n+1},B\power{u_i}) \\
     {\scriptstyle \Upsilon(g)}
     \bigg\downarrow & & \bigg\downarrow 
     {\scriptstyle \Upsilon(g)} \\
     (\widetilde{F}_{n+1},B\power{u_i}) &
     \stackrel{(X,\theta(h))}{\longrightarrow} &
     (\widetilde{F}_{n+1},B\power{u_i}). \\ 
    \end{array}\]
Note that the induced diagram on the residue field $L$
commutes.

Since $u_n^g\in E_{n+1}^0\subset B\power{u_i}$,
$(\theta(h)\circ \upsilon(g))(u_n)=u_n^g=
(\upsilon(g)\circ \theta(h))(u_n)$.
Hence $(\theta(h)\circ \upsilon(g))|_A=
(\upsilon(g)\circ \theta(h))|_A$.
From the fact that $B(i)$ is an \'{e}tale $A$-algebra,
$B\power{u_i}$ is complete, and the induced homomorphisms
on the residue field coincide,
we see that $(\theta(h)\circ \upsilon(g))|_{B(i)}=
(\upsilon(g)\circ \theta(h))|_{B(i)}$ for all $i$.
Hence $(\theta(h)\circ \upsilon(g))|_{B({\infty})}=
(\upsilon(g)\circ \theta(h))|_{B({\infty})}$ and 
$(\theta(h)\circ \upsilon(g))|_B=
(\upsilon(g)\circ \theta(h))|_B$.
Then the corollary follows from the universality of
$(\widetilde{F}_{n+1},B\power{u_i})$.
\qqq

\begin{theorem}\label{thm:Characteristic-zero}
The profinite group $\mathbb{G}$ acts on
the multiplicative cohomology theories 
$\mathbb{C}^*(-)$ and $\mathbb{B}^*(-)$.
There is an natural isomorphism of 
multiplicative cohomology theories
\[ {\mathbb B}^*(X)\cong \mathbb{C}^*(X),\]
with $\bigG$-action for all spectra $X$.
\end{theorem}

\proof
This follows from Proposition~\ref{operation-FGL}
and Corollary~\ref{cor:commutativity-actions}.
\qqq

There is a ring spectrum map
$E_{n+1}\to L_{K(n)}E_{n+1}=\mathbb{A}\to \mathbb{B}$.
By Lemma~\ref{lemma:description-C-theory} and
Theorem~\ref{thm:Characteristic-zero},
this induces a multiplicative natural transformation
of cohomology theories
\begin{align}\label{Chern-character} 
   \Theta :E_{n+1}^*(-)\longrightarrow 
           {\mathbb B}^*(-)
          \cong\mathbb{C}^*\widehat{\otimes}_{E_n^*}E_n^*(-). 
\end{align}
Note that $\Theta$ is equivariant with respect to
the action of $\mathbb{G}$ when we consider that 
$\mathbb{G}$ acts on the left hand side through
the projection $\mathbb{G}\to G_{n+1}$.
We say that $\Theta$ is a generalized Chern character
since it is a multiplicative natural transformation
from the height $(n+1)$ chomology theory $E_{n+1}^*(-)$
to the height $n$ cohomology theory 
$\mathbb{B}^*(-)\cong\mathbb{C}^*(-)$,
which is a coefficient extension of $E_n^*(-)$. 
Furthermore,
the inclusion $E_n^*\hookrightarrow \mathbb{C}^*$ induces 
a multiplicative natural transformation of 
cohomology theories
\[ I: E_n^*(-)\longrightarrow 
      \mathbb{C}^*\widehat{\otimes}_{E_n^*}E_n^*(-)
      \cong \mathbb{B}^*(-). \]
Then $I$ is equivariant with respect to
the action of $\mathbb{G}$ when we consider that 
$\mathbb{G}$ acts on the left hand side through
the projection $\mathbb{G}\to G_n$.

\section{The cohomology group $H_c^*(S_n;\mathbb{B}^*(X))$}
\label{section:comparison:Galois-extension}

In this section we show that the continuous cohomology
of the $n$th Morava stabilizer group $S_n$ with coefficients
in the cohomology group $\mathbb{B}^*(X)$
is naturally isomorphic to $\mathbb{A}^*(X)$
for every spectrum $X$.

In this section 
we give the graded commutative rings
$\mathbb{A}^*$ and $\mathbb{B}^*$ the $I_n$-adic topology.
In particular, the degree $0$ residue fields $K$ and $L$ 
are discrete.
Let $h=\mathbb{A}$ or $\mathbb{B}$.
Since $h^0$ is a complete Noetherian local ring, 
it is linealy compact.
By Lemma~\ref{lemma:inverselimit-description},
$h^*(X)\cong\ \subrel{\longleftarrow}{\lim}\,
     h^*(X_{\lambda})$,
where the inverse limit on the right hand side
is taken over $\Lambda(X)$. 
Let $\{F^{\delta}h^*(X)\}_{\delta\in \Delta(X)}$ be a 
family of $h^*$-submodules of $h^*(X)$ consisting
of the kernels of $h^*(X)\to h^*(Z)$,
which is induced by some map from a finite spectrum
$Z$ to $X$.
Then $\Delta(X)$ is a directed set by the reverse order of inclusions.
For any $(X_{\lambda}\to X)\in\Lambda(X)$,
there is a unique $\delta\in\Delta(X)$
such that $F^{\delta}h^*(X)$ is the kernel of
$h^*(X)\to h^*(X_{\lambda})$.
This defines a functor $\Lambda(X)\to \Delta(X)$;
$\lambda\mapsto\delta(\lambda)$,
which is final
(cf. \cite[IX.3]{MacLane}).
Hence we see that
$h^*(X) \cong\ \subrel{\longleftarrow\delta}{\lim}\,
     h^*(X)/F^{\delta}h^*(X)$.
In this section
we regard $h^*(X)$ as a filtered module by the filtration
\begin{align}\label{eq:filtration-h(X)} 
    \{ F^{\delta}h^*(X)+I_n^rh^*(X)
    \}_{(\delta,r)\in \Delta(X)\times\mathbf{N}_0}, 
\end{align}
where $\mathbf{N}_0$ is the set of non-negative integers.
Note that $h^*(X)/(F^{\delta}h^*(X)+I_n^rh^*(X))$
is a finitely generated Artinian $h^*$-module.
Then $h^*(X)/F^{\delta}h^*(X)$ is isomorphic to
$\subrel{\longleftarrow r}{\lim}\,h^*(X)/
   (F^{\delta}h^*(X)+I_n^rh^*(X))$.
This implies that 
$h^*(X)$ 
is isomorphic to 
$    \subrel{\longleftarrow\delta}{\lim}\,
     \subrel{\longleftarrow r}{\lim}\,
     h^*(X)/(F^{\delta}h^*(X)+I_n^rh^*(X))$.
Hence the filtration~(\ref{eq:filtration-h(X)}) 
of $h^*(X)$ is complete Hausdorff.

%
%
%

Let $\mathbb{B}(j)^0=B(j)\power{u_i}$ for $j\ge -1$.
By Lemma~\ref{henselian},
the sequence of Galois extensions~(\ref{sequence-Galois-extesion-residue})
induces a sequence of finite Galois extensions of commutative rings:
\[ \mathbb{A}^0=\mathbb{B}({-1})^0\to
              \mathbb{B}(0)^0\to
              \mathbb{B}(1)^0\to\cdots,\]
and the Galois group for the extension
$\mathbb{A}^0\to\mathbb{B}(i)^0$
is $S_n(i)$.
Let $\mathbb{B}({\infty})^0$ be
the direct limit of the sequence:
$\mathbb{B}({\infty})^0=\ \subrel{\longrightarrow i}{\lim}\,
\mathbb{B}(i)^0=\cup_i
\mathbb{B}(i)^0.$
Then $\mathbb{B}^0$ is isomorphic to be the $I_n$-adic completion
of $\mathbb{B}({\infty})^0$: 
\[ \mathbb{B}^0\cong (\mathbb{B}({\infty})^0)^{\wedge}_{I_n}. \]

Recall that an $\mathbb{A}^*$-module $M$ is pro-free if
$M$ is the $I_n$-adic completion for some free module. 
Let $S=\{s=(s_0,s_1,\ldots)\}$ be a set of multi-indexes such that 
$0\le s_0< p^n, 0\le s_i\le p^n\ (i\ge 1)$.
For $s\in S$, we set 
\[ \widetilde{\Phi}(s)=
   \widetilde{\Phi}_0^{s_0}\widetilde{\Phi}_1^{s_1}\cdots.\]
Then $\mathbb{B}({\infty})^*$ is a free $\mathbb{A}^*$-module
with basis $\{\widetilde{\Phi}(s)\}_{s\in S}$.
Hence we obtain the following lemma.

\begin{lemma}
$\mathbb{B}^*$ is pro-free over $\mathbb{A}^*$
with topological basis 
$\{\widetilde{\Phi}(s)\}_{s\in S}$.
\end{lemma}

\begin{lemma}\label{lemma:inverselimit-cokernel}
Let $\Lambda$ be an essentially small filtered category,
and let $\{M_{\lambda}\}_{\lambda\in\Lambda}$ be 
a $\Lambda$-diagram of finitely generated $\mathbb{A}^*$-modules
such that $\subrel{\longleftarrow}{\lim}M_{\lambda}=0$.
Then $\subrel{\longleftarrow}{\lim}
(\mathbb{B}^*\otimes_{\mathbb{A}^*}M_{\lambda})=0$.
\end{lemma}

\proof
Since $\mathbb{B}^*\otimes_{\mathbb{A}^*}M_{\lambda}$
is a finitely generated $\mathbb{B}^*$-module,
$\mathbb{B}^*\otimes_{\mathbb{A}^*}M_{\lambda}
\cong \mathbb{B}^*\widehat{\otimes}_{\mathbb{A}^*}M_{\lambda}
\cong \mathbb{B}({\infty})^*\widehat{\otimes}_{\mathbb{A}^*}M_{\lambda}$.
Then there is an injection
\[ \mathbb{B}({\infty})^*\widehat{\otimes}_{\mathbb{A}^*}M_{\lambda}
     \cong
    \left( \bigoplus_S M_{\lambda}\{\widetilde{\Phi}(s)\}
     \right)^{\wedge}_{I_n}
   \longrightarrow
    \prod_S M_{\lambda}.\]
This induces an injection
\[ \inverselimit{\lambda}(\mathbb{B}^*\otimes_{\mathbb{A}^*}
      M_{\lambda})\longrightarrow
   \inverselimit{\lambda}\left(\prod_S M_{\lambda}\right)
   \cong \prod_S \left(\inverselimit{\lambda}M_{\lambda}\right)
    =0.\]
Hence $\inverselimit{\lambda}(\mathbb{B}^*\otimes_{\mathbb{A}^*}
      M_{\lambda})=0$.
\qqq

\begin{lemma}\label{lemma:completion}
$\mathbb{B}^*/I_n^r\mathbb{B}^*\cong 
\mathbb{B}({\infty})^*/I_n^r\mathbb{B}({\infty})^*$ for every $r$.
\end{lemma}

\proof
See the proof of 
\cite[Theorem~A.1]{Hovey-Strickland}. 
\qqq

\begin{proposition}\label{prop:A-B-exchange}
There is a natural isomoprhism
$\mathbb{B}^*(X)\cong \mathbb{B}^*
\widehat{\otimes}_{\mathbb{A}^*}\mathbb{A}^*(X)$
for every spectrum $X$.
\end{proposition}

\proof
If $X$ is finite,
then $\mathbb{B}^*\otimes_{\mathbb{A}^*}\mathbb{A}^*(X)$
is finitely generated over $\mathbb{B}^*$.
Hence $\mathbb{B}^*\widehat{\otimes}_{\mathbb{A}^*}\mathbb{A}^*(X)
\cong \mathbb{B}^*\otimes_{\mathbb{A}^*}\mathbb{A}^*(X)\cong
\mathbb{B}^*(X)$.
For general $X$, by definition,
\[  \mathbb{B}^*\widehat{\otimes}_{\mathbb{A}^*}
    \mathbb{A}^*(X)\cong
    \subrel{\longleftarrow \delta,r}{\lim}(\mathbb{B}^*/I_n^r\mathbb{B}^*)
     \otimes_{\mathbb{A}^*}(\mathbb{A}^*(X)/F^{\delta}\mathbb{A}^*(X)). \]
Since $\mathbb{A}^*(X)/F^{\delta}\mathbb{A}^*(X)$
is finitely generated over $\mathbb{A}^*$,
$\subrel{\longleftarrow r}{\lim}(\mathbb{B}^*/I_n^r\mathbb{B}^*)
   \otimes_{\mathbb{A}^*}(\mathbb{A}^*(X)/F^{\delta}\mathbb{A}^*(X))$
is isomorphic to
$\mathbb{B}^*
   \otimes_{\mathbb{A}^*}
   (\mathbb{A}^*(X)/F^{\delta}\mathbb{A}^*(X))$.
By the fact that $\Lambda(X) \to \Delta(X)$ is a final functor,
\[ \subrel{\longleftarrow\delta}{\lim}\mathbb{B}^*
   \otimes_{\mathbb{A}^*}(\mathbb{A}^*(X)/F^{\delta}\mathbb{A}^*(X))
   \cong
   \subrel{\longleftarrow\lambda}{\lim}\mathbb{B}^*
   \otimes_{\mathbb{A}^*}(\mathbb{A}^*(X)/F^{\delta(\lambda)}
   \mathbb{A}^*(X)).\]
Let $M_{\lambda}$ be the cokernel of
the injection $\mathbb{A}^*(X)/F^{\delta(\lambda)}\mathbb{A}^*(X)
\to \mathbb{A}^*(X_{\lambda})$.
Then $\subrel{\longleftarrow}{\lim}M_{\lambda}=0$
since $\subrel{\longleftarrow}{\lim}
\mathbb{A}^*(X)/F^{\delta(\lambda)}\mathbb{A}^*(X)
\stackrel{\cong}{\to}\
\subrel{\longleftarrow}{\lim}
\mathbb{A}^*(X_{\lambda})$.
By Lemma~\ref{lemma:inverselimit-cokernel},
$\subrel{\longleftarrow}{\lim}
\mathbb{B}^*\otimes_{\mathbb{A}^*}M_{\lambda}=0$.
Since $\mathbb{A}^*\to\mathbb{B}^*$ is flat
by \cite[Lemma~23.1]{Matsumura},
this implies that 
$\subrel{\longleftarrow\lambda}{\lim}\mathbb{B}^*
 \otimes_{\mathbb{A}^*}(\mathbb{A}^*(X)/F^{\delta(\lambda)}
 \mathbb{A}^*(X))\cong \mathbb{B}^*(X)$
\qqq

Let $N$ be a $\mathbb{B}^*$-module with
filtration $\{F^{\delta}N\}_{\delta\in \Delta}$,
where $\Delta$ is a directed set given 
by the reverse order of inclusion.
For each $\delta$
we suppose that 
$N/F^{\delta}N$ is finitely generated and 
$I_n^rN/F^{\delta}N=0$ for some $r$. 
Furthermore,
we suppose that the filtration is complete Hausdorff.
Let $\mbox{\rm Map}_c(S_n,N)$ be the set of all
continuous maps from $S_n$ to the filtered module $N$.
Then $\mbox{\rm Map}_c(S_n,N)$ is a $\mathbb{B}^*$-module
by using the $\mathbb{B}^*$-module structure on $N$. 
Note that 
there is an exact sequence of $\mathbb{B}^*$-modules   
\[ 0\to\mbox{\rm Map}_c(S_n,F^{\delta}N)\longrightarrow
       \mbox{\rm Map}_c(S_n,N)\longrightarrow
       \mbox{\rm Map}_c(S_n,N/F^{\delta}N)\to 0\]
for any $\delta\in \Delta$.
We define a filtration on $\mbox{\rm Map}_c(S_n,N)$
by 
$\{\mbox{\rm Map}_c(S_n,F^{\delta}N)\}_{\delta\in
       \Delta}$.
Since there is an isomorphism
$\mbox{\rm Map}_c(S_n,N)\stackrel{\cong}{\to}
   \subrel{\longleftarrow \delta}{\lim} 
   \mbox{\rm Map}_c(S_n,N/F^{\delta}N)$,
this filtration is complete Hausdorff.
For $n\otimes b\in N{\otimes}_{\mathbb{A}^*}\mathbb{B}^*$,
the map $g\mapsto n\cdot g(b)$ 
defines a map $N\otimes_{\mathbb{A}^*}\mathbb{B}^*\to
\mbox{\rm Map}(S_n,N)$.
By Lemma~\ref{lemma:completion},
$b$ mod $I_n^r\mathbb{B}^*$ is contained in
$\mathbb{B}(i)^*/I_n^r\mathbb{B}(i)^*$ for some $i$.
Then $g(b)\equiv b$ mod $I_n^r\mathbb{B}^*$ for $g\in S_n^{(i)}$.
This shows that the map
$N\otimes_{\mathbb{A}^*}\mathbb{B}^*\to
\mbox{\rm Map}(S_n,N)$ factors through
$\mbox{\rm Map}_c(S_n,N)$.  
Then this map extends to a map
\[ N\widehat{\otimes}_{\mathbb{A}^*}\mathbb{B}^*
   \longrightarrow \mbox{\rm Map}_c(S_n,N). \]
Note that the map is a $\mathbb{B}^*$-module homomorphism
if we regard the left hand side as a $\mathbb{B}^*$-module
by using the $\mathbb{B}^*$-module structure on $N$.

\begin{proposition}\label{prop:mapping-exchange}
$N\widehat{\otimes}_{\mathbb{A}^*}\mathbb{B}^*
   \stackrel{\cong}{\longrightarrow} \mbox{\rm Map}_c(S_n,N)$.
\end{proposition}

\proof
If $N=\mathbb{B}^*$ with filtration
$\{I_n^r\mathbb{B}^*\}_{r\ge 0}$,
then by Lemma~\ref{lemma:completion}
\[ \begin{array}{rcl}
     \mathbb{B}^*\widehat{\otimes}_{\mathbb{A}^*}
     \mathbb{B}^*
     &\cong&
     \inverselimit{r}(\mathbb{B}^*/I_n^r\mathbb{B}^*)
     \otimes_{\mathbb{A}^*}(\mathbb{B}^*/I_n^r\mathbb{B}^*)\\
     &\cong&
     \inverselimit{r}
     (\mathbb{B}({\infty})^*/I_n^r\mathbb{B}({\infty})^*)
     \otimes_{\mathbb{A}^*}
     (\mathbb{B}({\infty})^*/I_n^r\mathbb{B}({\infty})^*)\\   
     &\cong&
     \inverselimit{r}\directlimit{i}
     (\mathbb{B}(i)^*/I_n^r\mathbb{B}(i)^*)
     \otimes_{\mathbb{A}^*}
     (\mathbb{B}(i)^*/I_n^r\mathbb{B}(i)^*).\\   
   \end{array}\] 
Since $\mathbb{A}^*\to\mathbb{B}(i)^*$ is a Galois extension
with Galois group $S_n(i)$ 
there is an isomorphism
\[ (\mathbb{B}(i)^*/I_n^r\mathbb{B}(i)^*)
   \otimes_{\mathbb{A}^*}
   (\mathbb{B}(i)^*/I_n^r\mathbb{B}(i)^*) \cong
   \mbox{\rm Map}(S_n(i),
   \mathbb{B}(i)^*/I_n^r\mathbb{B}(i)^*),\]
which is given by $x\otimes y\mapsto (g\mapsto x\cdot g(y))$.
Hence we obtain
\[ \begin{array}{rcl}
     \mathbb{B}^*\widehat{\otimes}_{\mathbb{A}^*}
     \mathbb{B}^*
     &\cong&
     \inverselimit{r}\directlimit{i}
     \mbox{\rm Map}(S_n(i),
     \mathbb{B}(i)^*/I_n^r\mathbb{B}(i)^*)\\
     &\cong&
     \inverselimit{r}
     \mbox{\rm Map}_c(S_n,
     \mathbb{B}({\infty})^*/I_n^r\mathbb{B}({\infty})^*)\\
     &\cong&
     \mbox{\rm Map}_c(S_n, \mathbb{B}^*).\\
   \end{array}\]
For $N$ general,
we have a sequence of isomorphisms: 
\[  \begin{array}{rcl} 
      N\widehat{\otimes}_{\mathbb{A}^*}\mathbb{B}^*
      &\cong &
      N\widehat{\otimes}_{\mathbb{B}^*}\mathbb{B}^*
      \widehat{\otimes}_{\mathbb{A}^*}\mathbb{B}^*\\[2mm]
      &\cong&
      N\widehat{\otimes}_{\mathbb{B}^*}
      \mbox{\rm Map}_c(S_n,\mathbb{B}^*)
      \\[2mm]
      &\cong&
     \mbox{\rm Map}_c(S_n,N).\\     
   \end{array}\] 
This completes the proof.
\qqq

We set 
$\mathbb{B}^{*\otimes p}
    =\overbrace{\mathbb{B}^*\widehat{\otimes}_{\mathbb{A}^*}\cdots 
    \widehat{\otimes}_{\mathbb{A}^*}\mathbb{B}^*}^p$
and $S_n^p=\overbrace{S_n\times \cdots\times S_n}^p$.
We define a continuous $\mathbb{B}^*$-algebra map
$\mathbb{B}^{*\otimes (p+1)}\to \mbox{\rm Map}_c(S_n^p,
\mathbb{B}^*)$
by 
\[ (b_0\otimes b_1\otimes \cdots \otimes b_p)\mapsto
 ((g_1,\ldots,g_p)\mapsto \prod_{i=0}^pg_1\cdots g_i(b_i)), \]
where the $\mathbb{B}^*$-algebra structure on 
$\mathbb{B}^{*\otimes (p+1)}$ comes from the first left factor.
By applying Proposition~\ref{prop:mapping-exchange} 
repeatedly, we obtain the following corollary.

\begin{corollary}\label{cor:iterated-exchange}
$\mathbb{B}^{*\otimes (p+1)}\stackrel{\cong}{\to}
\mbox{\rm Map}_c(S_n^p, \mathbb{B}^*)$
is an isomorphism.
\end{corollary}


\begin{corollary}\label{cor:general-map-exchange}
There is an isomorphism 
$\mathbb{B}^{*\otimes p} \widehat{\otimes}_{\mathbb{A}^*}N
\stackrel{\cong}{\to} \mbox{\rm Map}_c(S_n^p,N)$,
which is given by
\[ b_0\otimes \cdots \otimes b_{p-1}\otimes n \mapsto
    \left((g_1,\ldots,g_p)\mapsto 
    \left(\prod_{i=0}^{p-1}g_1\cdots g_i(b_i)\right)
    \cdot g_1\cdots g_p(n)\right). \]
\end{corollary}

For a topological module $M$ with
continuous $S_n$-action,
we denote by $C^p_c(S_n; M)$ the set of all continuous maps from
$S_n^p$ to $M$:
\[ C^p_c(S_n; M)=
   \mbox{\rm Map}_c(S_n^p,M). \]
Then we define a differential 
$d^p: C^p_c(S_n; M)\longrightarrow C^{p+1}_c(S_n; M)$
by
\[ \begin{array}{rcl} 
     d^p(\varphi)(g_1,\ldots,g_{p+1}) & = &
       g_1\cdot\varphi(g_2,\ldots,g_{p+1})  \\[2mm]
    &&\quad + \sum_{i=1}^p 
      (-1)^i\varphi(g_1,\ldots,g_ig_{i+1},\ldots,g_{p+1})\\[2mm]
    &&\quad +(-1)^{p+1} \varphi(g_1,\ldots,g_p).
   \end{array} \]
Then $C^*_c(S_n;M)=\{C^p_c(S_n;M),d^p\}_{p\ge 0}$
forms a cochain complex and the continuous cohomology of
$S_n$ with coefficients in $M$ is defined to be the 
cohomology of $C^*_c(S_n;M)$: 
\[ H^p_c(S_n;M):= H^p(C^*_c(S_n;M)) .\]

Since the Morava stabilizer group $S_n$ 
acts on $\mathbb{B}^*(X)$ as cohomology operations,
$\mathbb{B}^*(X)$ is a complete Hausdorff 
filtered $\mathbb{B}^*$-module with compatible 
continuous $S_n$-action.

\begin{proposition}\label{prop:fundamental_exchange}
There is a natural isomorphism 
of cochain complexes:
\[ C_c^*(S_n; \mathbb{B}^*(X))\cong 
   C_c^*(S_n, \mathbb{B^*})\widehat{\otimes}_{\mathbb{A}^*}
   \mathbb{A}^*(X)\]
for every spectrum $X$.
\end{proposition}

\proof
By Proposition~\ref{prop:A-B-exchange},
there is a natural isomorphism
$\mathbb{B}^*(X)\cong\mathbb{B}^*
\widehat{\otimes}_{\mathbb{A}^*}\mathbb{A}^*(X)$.
This isomorphism is actually an isomorphism
of $S_n$-modules, 
where the $S_n$-action on $\mathbb{A}^*(X)$ is trivial.
Then the proposition
follows from Corollary~\ref{cor:iterated-exchange} and 
Corollary~\ref{cor:general-map-exchange}.
\qqq

By taking the coefficient of 
$1=\widetilde{\Phi}_0^0\widetilde{\Phi}_1^0\cdots$
with respect to the topological basis 
$\{\widetilde{\Phi}(s)\}_{s\in S}$,
there is a continuous $\mathbb{A}^*$-module homomorphism
$\varepsilon: \mathbb{B}^*\to \mathbb{A}^*$,
which gives a splitting of the unit 
$\eta: \mathbb{A}^*\to \mathbb{B}^*$.  
By Corollary~\ref{cor:iterated-exchange},
$C_c^p(S_n;\mathbb{B}^*)=\mbox{\rm Map}_c(S_n^p,\mathbb{B}^*)
\cong \mathbb{B}^{*\otimes (p+1)}$.
Then $d^p$ corresponds to the map
$d^p: \mathbb{B}^{*\otimes (p+1)}\to \mathbb{B}^{*\otimes (p+2)}$
given by
\[ d^p(b_0\otimes b_1\otimes \cdots \otimes b_p)=
   \sum_{i=0}^{p+1}(-1)^i b_0\otimes\cdots b_{i-1}\otimes 1
   \otimes b_i\otimes \cdots b_p.\]
Define a continuous $\mathbb{A}^*$-module map
$s^p: \mathbb{B}^{*\otimes (p+2)}\to \mathbb{B}^{*\otimes (p+1)}$
by 
\[ s^p(b_0\otimes\cdots \otimes b_{p+1})= 
   \varepsilon (b_0)\cdot b_1\otimes \cdots \otimes b_{p+1}, \]
where $\varepsilon: \mathbb{B}^*\to \mathbb{A}^*$
is the splitting of the inclusion $\eta: \mathbb{A}^*\to \mathbb{B}^*$.

\begin{lemma}\label{lemma:split_exact_sequence}
The sequence  
\begin{equation}\label{lema:split-exact-fundamental}
 0\to \mathbb{A}^* \stackrel{\eta}{\longrightarrow} \mathbb{B}^*
   \stackrel{d^0}{\longrightarrow} \mathbb{B}^{* \otimes 2}
   \stackrel{d^1}{\longrightarrow} \mathbb{B}^{*\otimes 3}
   \longrightarrow \cdots 
\end{equation}
is a split exact sequence of 
complete Hausdorff filtered $\mathbb{A}^*$-modules.
\end{lemma}

\proof
It is easy to check that
$s^pd^p+d^{p-1}s^{p-1}$
is the identity map of $\mathbb{B}^{*\otimes (p+1)}$ for $p>0$.
Furthermore,
$s^0d^0+\eta\varepsilon=id_{\mathbb{B}^*}$ and
$\varepsilon\eta=id_{\mathbb{A}^*}$.  
These show that 
(\ref{lema:split-exact-fundamental})
is a split exact sequence.
\qqq

By applying the functor 
$(-)\widehat{\otimes}_{\mathbb{A}^*}\mathbb{A}^*(X)$
to the split exact sequence (\ref{lema:split-exact-fundamental})
and by using Proposition~\ref{prop:fundamental_exchange},
we obtain the following corollary. 

\begin{corollary}
For every spectrum $X$, the sequence 
\[ 0\to \mathbb{A}^*(X)\stackrel{\eta}{\longrightarrow}
        C^0_c(S_n;\mathbb{B}^*(X))\stackrel{d^0}{\longrightarrow}
        C^1_c(S_n;\mathbb{B}^*(X))\stackrel{d^1}{\longrightarrow}         
        C^2_c(S_n;\mathbb{B}^*(X))\stackrel{d^2}{\longrightarrow}
         \cdots \]
is a split exact sequence of complete Hausdorff
filtered $\mathbb{A}^*$-modules.
\end{corollary}

Hence we obtain the following theorem.

\begin{theorem}\label{theorem:S_n-cohomology}
For every spectrum $X$, 
$H^i_c(S_n; \mathbb{B}^*(X))=0$ for $i>0$,
and there is a natural isomorphism
of $G_{n+1}$-modules
\[ H^0_c(S_n; \mathbb{B}^*(X))\cong \mathbb{A}^*(X). \]
\end{theorem}

\section{Comparison of $\mathbb{B}^*(X)$ and $E_n^*(X)$}
\label{section:comparison->En}

By Lemma~\ref{lemma:Gn-module-inclusion},
there is a natural injective map
of $G_n$-modules
$E_n^*(X)\longrightarrow H^0(S_{n+1}; \mathbb{B}^*(X))$.
In this section we prove the folloing theorem.

\begin{theorem}[Theorem~A]\label{thm:invariantS_n+1}
For every spectrum $X$,
there is a natural isomorphism of
$E_n^*$-modules with compatible continuous  
$G_n$-action
\[ E_n^*(X)\stackrel{\cong}{\longrightarrow} 
              H^0(S_{n+1}; \mathbb{B}^*(X)). \]
For a finite spectrum $X$,
$E_n^*(X)$ with $G_n$-action can be recovered from
$E_{n+1}^*(X)$ with $G_{n+1}$-action.
\end{theorem}

Let $S$ be a complete Noetherian local ring,
and $R$ a subring which is also a complete Noetherian
local ring.
Suppose that the inclusion $R\hookrightarrow S$
is a flat local homomorphism,
and the maximal ideal $\mathfrak{m}_S$ of $S$ is 
generated by the maximal ideal $\mathfrak{m}_R$ of $R$:
$\mathfrak{m}_R S = \mathfrak{m}_S$.
In this case $S$ is faithfully flat over $R$
by \cite[Theorem~7.2]{Matsumura}.
Hence $R/\mathfrak{m}_R^i\to S/\mathfrak{m}_R^iS=S/\mathfrak{m}_S^i$
is also faithfully flat by base change,
and in particular injective.   
Let $G$ be a group acting on $S$ as ring automorphisms.
Suppose that $R$ is contained in the fixed subring $S^G=H^0(G;S)$. 

\begin{lemma}\label{lemma:R-S-invariant}
If $R/\mathfrak{m}_R=(S/\mathfrak{m}_S)^G$,
then $R=S^G$.
\end{lemma}

\proof
Since $H^0(G;S)\cong\ \subrel{\longleftarrow}{\lim}
H^0(G;S/\mathfrak{m}_S^i)$,
it is sufficient to show that $R/\mathfrak{m}_R^i=(S/\mathfrak{m}_S^i)^G$
for all $i$.
We prove this by induction on $i$.
Assume that the $i-1$ case is true.
There is a commutative diagram:
\[ \begin{array}{ccccccccc}
     0&\to& \mathfrak{m}_R^{i-1}/\mathfrak{m}_R^i &
     \longrightarrow & R/\mathfrak{m}_R^i &
     \longrightarrow & R/\mathfrak{m}_R^{i-1} & \to &0 \\
     && \bigg\downarrow & & \bigg\downarrow && \bigg\downarrow && \\ 
     0&\to& (\mathfrak{m}_S^{i-1}/\mathfrak{m}_S^i)^G &
     \longrightarrow & (S/\mathfrak{m}_S^i)^G &
     \longrightarrow & (S/\mathfrak{m}_S^{i-1})^G. & & \\
   \end{array}\]
Note that the middle vertical map is injective,
and the right vertical map is an isomorphism
by the hypothesis of induction.
Since $\mathfrak{m}_S=\mathfrak{m}_R S$,
the $S/\mathfrak{m}_S$-module 
$\mathfrak{m}_S^{i-1}/\mathfrak{m}_S^i$ is generated by
the $(i-1)$ products of generators of $\mathfrak{m}_R$.
Hence we can take a basis of 
$\mathfrak{m}_S^{i-1}/\mathfrak{m}_S^i$
consisting of elments in the image of the left vertical map.
Then $(S/\mathfrak{m}_S)^G=R/\mathfrak{m}_R$
implies that the left vertical map is surjctive.
Hence we see that the middle vertical is an isomorphism.
\qqq

\begin{lemma}\label{lemma:invariant-basic}
$H^0(S_{n+1};L[u^{\pm 1}])=\mathbf{F}[w^{\pm 1}]$.
\end{lemma}

\proof
Let $M^1_n=v_{n+1}^{-1}BP_*/(p,v_1,\ldots,v_{n-1},v_n^{\infty})$.
By \cite[Theorem~5.10]{MRW},
$\mbox{\rm Ext}^0_{BP_*(BP)}(BP_*,M^1_n)$ is
the direct sum of the finite torsion submodules
and the $K(n)_*/k(n)_*$ generated by $1/v_n^j,\ j\ge 1$ 
as a $k(n)_*$-module.
Then as in \cite[\S5.3]{Torii1}
$H^0(S_{n+1};\mathbf{F}\power{u_n}[u^{\pm 1}])={\mathbf F}[v_n]$.
By \cite[Lemma~5.9]{Torii1},
$H^0(S_{n+1};{\mathbf F}((u_n))[u^{\pm 1}])$ is the localization of
$H^0(S_{n+1};\mathbf{F}\power{u_n}[u^{\pm 1}])$ 
by inverting the invariant element $v_n$.
Hence $H^0(S_{n+1};{\mathbf F}((u_n))[u^{\pm 1}])={\mathbf F}[v_n^{\pm 1}]$.
By \cite[Lemma~3.7]{Torii1}, 
$w=\Phi_0^{-1}u\in L$
is invariant under the action of $S_{n+1}$.
Let $a$ be a degree $2n$
invariant element in $L[u^{\pm 1}]$.
Then $b=a w^n\in L$
is also invariant.
Let $\phi(X)\in {\mathbf F}((u_n))[X]$ be the minimal polynomial of $b$.
Then $\phi(b)=0$.
Since $b$ is invariant under the action of $S_{n+1}$,
$\phi^g(b)=0$ for all $g\in S_{n+1}$.
Hence $\phi^g(X)$ is also the minimal polynomial of $b$.
This implies that $\phi(X)$ is a polynomial over
$H^0(S_{n+1};{\mathbf F}((u_n)))={\mathbf F}$.
Hence $b\in \overline{{\mathbf F}}\cap L={\mathbf F}$.
This completes the proof. 
\qqq

\begin{corollary}\label{cor:invariant-base-case}
For $0\le i\le n$,
$H^0(S_{n+1}; \mathbb{B}^*/I_i)\cong E_n^*/I_i$.
\end{corollary}

\proof
Since $E_n^*\to \mathbb{B}^*$ is a flat local homomorphism,
$E_n^*/I_i\to \mathbb{B}^*/I_i$ is faithfully flat.
In particular it is injective.
By Lemma~\ref{lemma:Gn-module-inclusion},
$E_n^* \subset H^0(S_{n+1}; \mathbb{B}^*)$.
This implies that
$E_n^*/I_i\subset H^0(S_{n+1}; \mathbb{B}^*/I_i)$.
Then the corollary follows from Lemma~\ref{lemma:R-S-invariant}
and Lemma~\ref{lemma:invariant-basic}.
\qqq

Let $M$ be a finitely generated $E_n^*$-module.
Note that $\mathbb{B}^*\widehat{\otimes}_{E_n^*}M
=\mathbb{B}^*\otimes_{E_n^*}M$.
The action of $S_{n+1}$ on $\mathbb{B}^*$ and the trivial action
on $M$ defines an action of $S_{n+1}$ on 
$\mathbb{B}^*\widehat{\otimes}_{E_n^*}M$.  

\begin{lemma}\label{lemma:invariant-sn+1induction}
Suppose that $M$ has a finite filtration
$M=M_0\supset M_1\supset \cdots\supset M_s=0$
such that each quotient $M_i/M_{i+1}$ is
isomorphic to (suspention of) $E_n^*/I_k$ for some $0\le k\le n$.
Then $H^0(S_{n+1}; \mathbb{B}^*\widehat{\otimes}_{E_n^*}M)=M$.
\end{lemma}

\proof
We prove the lemma by induction on the length of filtration.
Assume that it is true if the length of filtration is less than $s$. 
Suppose that $M$ has length $s$.
Then there is a exact sequence of $E_n^*$-modules: 
$0\to N\to M\to M/N\to 0$,
where $N$ has length $s-1$ and $M/N$ is isomorphic to
$E_n^*/I_k$.
This induces the following commutative diagram:
\[ \begin{array}{rcccccl}
    0 \to & N & \longrightarrow & M & \longrightarrow & M/N & \to 0\\
          & \bigg\downarrow & & \bigg\downarrow & & 
                                \bigg\downarrow &   \\
    0 \to & \phantom{1S_n}(\mathbb{B}^*\otimes N)^{S_{n+1}} & 
            \longrightarrow &
    \phantom{1S_n}(\mathbb{B}^*\otimes M)^{S_{n+1}} & \longrightarrow &
    \phantom{SN+1}(\mathbb{B}^*\otimes M/N)^{S_{n+1}}. & \\
   \end{array}\]
The right vertical arrow is an isomorphism by 
Lemma~\ref{cor:invariant-base-case},
and the left vertical arrow 
is an isomorphism by the hypothesis of induction. 
Hence the middle vertical arrow is also an isomorphism. 
\qqq

\begin{corollary}\label{cor:invariantS_n+1forfintie}
For any finite spectrum $X$,
$E_n^*(X)\stackrel{\cong}{\to}H^0(S_{n+1}; \mathbb{B}^*(X))$.
\end{corollary}

\proof
By Landweber filtration theorem,
there is a finite filtration of $BP^*$-modules:
$BP^*(X)=F_0\supset F_1\supset \cdots\supset F_s=0$
such that each quotient $F_i/F_{i+1}$
is isomorphic to (supsension of) $BP^*/I_k$ for some $k$.
Since $E_n^*(X)=E_n^*\otimes_{BP^*}BP^*(X)$,
there is a filtration 
$E_n^*(X)=F_0'\supset F_1'\supset\cdots\supset F_t'=0$ 
such that $F_i'/F_{i+1}'=E_n^*/I_k$ for some $0\le k\le n$.
Then the corollary follows from 
Lemma~\ref{lemma:invariant-sn+1induction}.
\qqq

\proof[Proof of Theorem~\ref{thm:invariantS_n+1}]
Since $E_n^*$ and $\mathbb{B}^*$ are complete Noetherian local rings,
$E_n^*(X)\cong\ \subrel{\longleftarrow}{\im}E_n^*(X_{\lambda})$
and 
$\mathbb{B}^*(X)\cong\ \subrel{\longleftarrow}{\lim}
\mathbb{B}^*(X_{\lambda})$.
Then Theorem~\ref{thm:invariantS_n+1} follows
from Corollary~\ref{cor:invariantS_n+1forfintie}
and the fact that 
$\mathbb{B}^*(X)\cong
\mathbb{B}^*\otimes_{E_{n+1}^*}E_{n+1}^*(X)$
for a finite spectrum $X$.
\qqq

\begin{remark}\rm
Let $p$ be an odd prime number.
By the same method as above,
we can show that 
there is a natural isomorphism
\[ (E_n/I_i)^*(X)\stackrel{\cong}{\longrightarrow}
   H^0(S_{n+1};(\mathbb{B}/I_i)^*(X))\]
of $G_{n}$-modules
for all spectra $X$ and $0\le i\le n$.
\end{remark}

For $m=n$ or $n+1$,
we define a pro-obeject $\mathbb{E}_m^*(X)$ 
of finitely generated $E_m^*$-modules with compatible 
continuous $G_m$-action
to be the system
\[ \mathbb{E}_m^*(X)=\{E_m^*(X_{\lambda})\}_{\lambda\in\Lambda(X)}. \]
Note that there is a natural isomorphism
$E_m^*(X)\stackrel{\cong}{\to}\
\subrel{\longleftarrow}{\lim}\mathbb{E}_m^*(X)$.
For a finite spectrum $X_{\lambda}$,
$\mathbb{B}^*(X_{\lambda})\cong \mathbb{B}^*\otimes_{E_{n+1}^*}
E_{n+1}^*(X_{\lambda})$.
By Theorem~\ref{thm:invariantS_n+1},
the fintiely generated $E_n^*$-module 
$E_n^*(X_{\lambda})$
with compatible continuous $G_n$-action
is isomorphic to $H^0(S_{n+1};
\mathbb{B}^*\otimes_{E_{n+1}^*}E_{n+1}^*(X_{\lambda}))$.
Hence we obtain the following corollary.

\begin{corollary}\label{cor:recover}
For any spectrum $X$,
the pro-object $\mathbb{E}_n^*(X)$
is isomorphic to the pro-object 
\[ H^0(S_{n+1};\mathbb{B}^*\otimes_{E_{n+1}^*}
   \mathbb{E}_{n+1}^*(X)) =
   \{H^0(S_{n+1};\mathbb{B}^*\otimes_{E_{n+1}^*}
   E_{n+1}^*(X_{\lambda}))\}. \]
\end{corollary}

Corollary~\ref{cor:recover} means that 
the pro-object $\mathbb{E}_n^*(X)$ can
be recovered from the pro-object $\mathbb{E}_{n+1}^*(X)$.

\begin{remark}\rm
The cohomology group $E_n^*(X)$ can not be obtained
from $E_{n+1}^*(X)$ in general,
since the cohomological Bousfield class $\langle E_n^*\rangle$
is $\langle K(n)\rangle$ (cf. \cite{Hovey2}).
For example,
$E_n^*(E_n)\neq 0$ but $E_{n+1}^*(E_n)=0$.
\end{remark}

{\footnotesize

}

\vspace{3mm}

\begin{flushleft}
Department of Mathematics, Faculty of Science,
Okayama University,\\
Okayama 700--8530, Japan\\
{\it E-mail}\,: torii@math.okayama-u.ac.jp\\
\end{flushleft}

\end{document}